\documentclass[12pt]{amsart}

\usepackage{amsmath,amssymb,amsfonts,amsthm,graphics,pst-plot,verbatim}

\newcommand{\nc}{\newcommand}
\nc{\dmo}{\DeclareMathOperator}
\nc{\nt}{\newtheorem}

\dmo{\cat}{CAT(0)}
\dmo{\teich}{T}
\dmo{\isom}{Isom}
\dmo{\aut}{Aut}
\dmo{\pants}{P}
\dmo{\cc}{C}
\dmo{\vis}{{\mathcal V}}
\dmo{\pd}{d_P}
\dmo{\dwp}{d_{WP}}
\dmo{\mcg}{Mod}
\dmo{\homeo}{Homeo}

\nc{\homeopms}{\homeo^{\pm}(S)}

\nc{\ep}{\epsilon}

\dmo{\gl}{GL}
\dmo{\pgl}{PGL}
\dmo{\slxxx}{SL}
\nc{\glz}{\gl_2({\mathbb Z})}
\nc{\pglz}{\pgl_2({\mathbb Z})}
\nc{\slz}{\slxxx_2({\mathbb Z})}

 \nc{\tbt}[4]{\ensuremath{\left(\begin{array}{rr} #1 & #2 \\
#3 & #4 \\ \end{array}\right)}}

\dmo{\f}{{\mathcal F}}
\nc{\ft}{\f}
\nc{\fkt}{\f_k}
\nc{\flt}{\f_l}
\nc{\fzt}{\f_0}
\nc{\fri}{\f_i}
\nc{\fgek}{\f_{\ge k}}
\nc{\fgk}{\f_{>k}}
\nc{\fit}{\f_i}
\nc{\fmk}{\f_M}
\nc{\fmok}{\f_{M-1}}
\nc{\fkmo}{\f_{k-1}}

\nc{\vxs}{\vis_X(S)}
\nc{\mods}{\mcg(S)}
\nc{\modf}{\mcg(R)}
\nc{\kb}{K'}
\nc{\h}{{\mathbb H}^2}
\nc{\coll}{{\mathcal C}}

\nc{\ib}{\bar{I}}
\nc{\ibs}{\ib_\star}

\nc{\ts}{\teich(S)}
\nc{\tsb}{\overline{\ts}}
\nc{\tsoo}{\teich(\soo)}
\nc{\tszf}{\teich(\szf)}

\nc{\cps}{\pants(S)}
\nc{\cpos}{\pants^0(S)}
\nc{\ccs}{\cc(S)}

\nc{\soo}{S_{1,1}}
\nc{\szf}{S_{0,4}}
\nc{\szfi}{S_{0,5}}
\nc{\sot}{S_{1,2}}
\nc{\stz}{S_{2,0}}

\nc{\tf}{\teich(R)}
\nc{\cpf}{\pants(R)}
\nc{\tfb}{\overline{\tf}}
\nc{\sooszfo}{R}

\nc{\fg}{{\mathcal G}}

\dmo{\os}{{\mathcal O}}

\nc{\set}[1]{\{#1\}}

\nc{\zt}{{\mathbb Z}_2}

\nc{\bl}{ \begin{list}{$\cdot$}{
\setlength{\leftmargin}{.25in}
\setlength{\rightmargin}{.5in}
\setlength{\parsep}{0.5ex plus .2ex minus 0ex}
\setlength{\itemsep}{0.2ex plus 0.2ex minus 0ex} 
}
}

\nc{\el}{\end{list}}

\nc{\pic}[2]{\begin{figure}[htb] \center{\leavevmode 
\epsfbox{#1.eps}} \caption{#2} \label{#1pic}\end{figure}} 

\nc{\pics}[3]{\epsfysize=#3 cm \begin{figure}[htb] 
\center{\leavevmode \epsfbox{#1.eps}} \caption{#2}  
\label{#1pic}\end{figure}} 

\nt{thm}{Theorem}[section]
\nt*{mthm}{Main Theorem}
\nt{prop}[thm]{Prop}
\nt{lem}[thm]{Lemma}
\nt{cor}[thm]{Corollary}
\nt*{q}{Question}

\nc{\bpf}{\begin{proof}}
\nc{\epf}{\end{proof}}

\nc{\p}[1]{{\bf #1}}

\include{diagram} 

\begin{document}

\input{epsf.sty}

\title{\bf Weil--Petersson isometries via the pants complex}

\author{Jeffrey Brock}

\author{Dan Margalit}

\address{Department of Mathematics\\ Box 1917\\Providence, RI 02912}

\address{Department of Mathematics\\ University of Utah\\ 155 S 1440 East \\ Salt Lake City, UT 84112-0090}

\thanks{The first author is partially supported by NSF grant number 0354288.  The second author is partially supported by an NSF postdoctoral fellowship and a VIGRE postdoctoral position under NSF grant number 0091675 to the University of Utah.}

\email{brock@math.brown.edu, margalit@math.utah.edu}

\keywords{Teichm\"uller space, Weil--Petersson metric, isometries, mapping class groups}

\subjclass[2000]{Primary: 32G15; Secondary: 32M99}

\maketitle

\begin{center}\today\end{center}

\begin{abstract}We extend a theorem of Masur and Wolf which says that given a hyperbolic surface $S$, every isometry of the Teichm\"uller space for $S$ with the Weil--Petersson metric is induced by an element of the mapping class group for $S$.  Our argument handles the previously untreated cases of the four-holed sphere, the one-holed torus, and the two-holed torus.\end{abstract}

\section{Introduction}

The purpose of this paper is to fill in several gaps in the understanding of Weil--Petersson isometries.
Throughout, $S$ is an orientable hyperbolic surface, and $S_{g,n}$ denotes a surface of genus $g$ with $n$ punctures.  The {\em extended mapping class group} for $S$ is: \[ \mods = \pi_0(\homeopms) \]
Let $\ts$ denote the Teichm\"uller space of $S$ with the Weil--Petersson metric.  Masur and Wolf proved that for $S \notin \set{\soo,\szf,\sot}$, every isometry of $\ts$ is induced by an element of $\mods$ \cite{mw}.  One of the main ingredients is that for exactly these surfaces, every automorphism of the complex of curves for $S$ is induced by $\mods$; this is a theorem of Ivanov \cite{nvi}, with special cases due to Korkmaz and Luo \cite{mk} \cite{fl}.  See Section~\ref{bg} for definitions.

We are able to extend the result of Masur and Wolf to all possible surfaces by viewing an isometry $I$ of $\ts$ as giving an automorphism $\ibs$ of the {\em pants graph} for $S$.  By a theorem of the second author, $\ibs$ is necessarily induced by a mapping class $f$.  Then, as observed by Wolpert, one can apply a theorem of the first author to show that $f$ induces $I$.  Thus, we have:

\begin{mthm}The natural map $\eta: \mods \to \isom(\ts)$ is surjective.  Further, $\ker(\eta) \cong \zt$ for $S \in \set{\soo,\sot,\stz}$, $\ker(\eta) \cong  \zt \oplus \zt$ for $S = \szf$, $\ker(\eta) = \mods$ for $S = S_{0,3}$ and $\ker(\eta)=1$ otherwise.\end{mthm}

It follows from the definition of the Weil--Petersson metric that $\mods$ acts by isometries; the content is that all isometries arise in this way.  

The application of curve complexes to the study of isometries of Teichm\"uller space first appears in Ivanov's proof of Royden's Theorem, which is the analogue of our Main Theorem for the Teichm\"uller metric \cite{nvi}.

\p{Acknowledgements.} We would like to thank Benson Farb for sharing many mathematical ideas; in particular, he suggested using the pants complex to understand Weil--Petersson isometries.  This paper came out of a minicourse hosted by the University of Utah; we would like to thank the participants, in particular Bob Bell and Ken Bromberg, for helpful discussions.  The remark at the end of Section~\ref{base} came from discussions with Chris Leininger.  Also, the proof of Lemma~\ref{strata} presented here is due to Mladen Bestvina; we thank him, as well as Howie Masur, for many enlightening conversations.  We are grateful to Indira Chatterji and Kevin Wortman for useful comments on earlier drafts of the paper.

\section{Preliminaries} \label{bg}

We briefly recall the basic definitions and concepts.

\p{Complex of curves.} The {\em complex of curves} $\ccs$ for $S$,
defined by Harvey \cite{wjh}, is the abstract simplicial flag complex
with vertices corresponding to isotopy classes of simple closed curves
in $S$, and edges between vertices that can be realized disjointly in
$S$.

For our purposes, this complex has an undesirable property when $S$ is $\soo$, $\szf$, or $\sot$: there are automorphisms of $\ccs$ which are not induced by elements of $\mods$.  For $S \in \set{\soo,\szf}$, this is because, as defined, $\ccs$ is a countable discrete set.  

The problem with $\sot$ is more subtle.  Luo noticed that if $\iota$
is the hyperelliptic involution of $\sot$, then the projection $\pi:
\sot \to \sot/\iota$ gives a bijection between the
vertices of $\cc(\sot)$ and $\cc(\sot/\iota) \cong \cc(\szfi)$.  Since
$\mcg(\szfi)$ acts transitively on the vertices of $\cc(\szfi)$ (every
curve has two punctures on one side and three on the other), it
follows that some elements of $\aut(\cc(\sot))$ fail to preserve the
set of vertices corresponding to separating curves.  Such
automorphisms clearly cannot arise from mapping classes (see
\cite{fl}).

\p{Pants graph.} The {\em pants graph} $\cps$ for $S$, introduced by Hatcher and Thurston \cite{ht} \cite{ah}, is the simplicial complex with vertices corresponding to {\em pants decompositions} of $S$ (i.e. maximal simplices of $\ccs$), and edges connecting pants decompositions which differ by an {\em elementary move}, by which we mean that the two pants decompositions differ by only one curve, and the differing curves have the smallest possible intersection.  The minimal intersection number is 1 or 2, depending on whether the curve being replaced lies in a single pair of pants or is the boundary between two pairs of pants; see Figure~\ref{movespic}.

\pics{moves}{Elementary moves.}{1.5}

We will apply the following theorem of the second author \cite{dm}:

\begin{thm}\label{dan}The natural map $\theta: \mods \to \aut(\cps)$ is surjective.  Further, $\ker(\theta) \cong \zt$ for $S \in \set{\soo,\sot,\stz}$, $\ker(\theta) \cong  \zt \oplus \zt$ for $S = \szf$, $\ker(\theta) = \mods$ for $S = S_{0,3}$, and $\ker(\theta)=1$ otherwise.\end{thm}

We note that the proof of Theorem~\ref{dan} is based on the aforementioned theorem of Ivanov, which gives an analogous statement for $\ccs$.

\p{Teichm\"uller space.} A point in the {\em Teichm\"uller space} $\ts$ of $S$ is given by a pair $(X,f)$, where $X$ is a finite-area hyperbolic surface $X$ and $f : S \to X$ is a homeomorphism.  Two points $(X,f)$ and $(Y,g)$ are equivalent if $g \circ f^{-1}$ is isotopic to an isometry.

We note that $\teich(S_{0,3})$ is a single point, and in this case our main result is trivial.

\p{Weil--Petersson metric.} A point in $\ts$ is naturally a Riemann surface via its uniformization as a quotient of $\h$ by a Fuchsian group.  In this way, the cotangent space $T_X^\star(\ts)$ at a point $X$ is identified with the space of holomorphic quadratic differentials on $X$ (holomorphic forms of type $\phi(z)dz^2$), which has the $L^2$ inner product defined by: \[ \langle \phi,\psi \rangle = \int_X \frac{\phi\bar{\psi}}{\rho^2} \] Here, $\rho(z)|dz|$ is the hyperbolic metric on $X$.  Then, the {\em Weil--Petersson metric} is defined by the pairing \[ \langle \mu,\phi \rangle = \int_X \mu \phi \] for $\mu \in \ts$ and $\phi \in T_X^\star(\ts)$.

Chu and Wolpert showed that the Weil--Petersson metric is not complete \cite{tc} \cite{sw1}.

\p{Augmented Teichm\"uller space.}  Masur gave an interpretation of points in the completion of $\ts$ as {\em marked noded Riemann surfaces} \cite{hm}.  A {\em noded Riemann surface} is a complex space with at most isolated singularities, called {\em nodes}, each possessing a neighborhood biholomorphic to a neighborhood of $(0,0)$ in the curve $\set{zw = 0}$ in ${\mathbb C}^2$.  Removing the nodes of a noded Riemann surface $W$ yields a (possibly disconnected) Riemann surface whose components we call the {\em pieces} of $W$.

Given a simplex $\sigma \in \ccs$, a {\em marked noded Riemann surface $(W,f)$ with nodes at $\sigma$} is a noded Riemann surface $W$ equipped with a continuous mapping $$f \colon S \to W$$ so that $f \vert_{S \setminus \sigma}$ is a homeomorphism to the pieces of $W$.  Two marked noded Riemann surfaces $(W,f)$ and $(Z,g)$ are equivalent if there is a continuous marking-preserving map $\phi \colon W \to Z$ which preserves nodes and is biholomorphic on the pieces.

\newcommand{\reals}{\mathbb{R}}

To describe a neighborhood of a point $(W,f)$ in $\tsb$, we give coordinates adapted to the simplex $\sigma$.  Given a maximal simplex $\tau$ with $\sigma \subset \tau$, there are \emph{Fenchel-Nielsen coordinates} 
\[ \{(\ell,\theta)_{\alpha}\} \in \left( \mathbb{R}_{>0} \times \mathbb{R}\right)^{|\tau|} \]
 for $\ts$ specifying length and twist parameters along each closed geodesic $\alpha$ determined by the vertices of $\tau$ (see, e.g. \cite{it}).  Then the \emph{extended Fenchel-Nielsen coordinates} for $\tau$ are obtained by allowing lengths to range in $\reals_{\ge 0}$ and stipulating that for any twist parameters $\theta$ and $\theta'$ we have $(0,\theta)_\alpha \sim (0,\theta')_\alpha.$ 

The extended Fenchel-Nielsen coordinates for $(W,f)$, then, are given by setting the length parameters corresponding to the curves of $\sigma$ to $0$, and setting the other length and twist parameters equal to the corresponding Fenchel-Nielsen coordinates for the piece in which the geodesic lies.  Then a neighborhood of $(W,f)$ is given by all surfaces in $\tsb$ with extended Fenchel-Nielsen coordinates close to those of $(W,f)$.  We will in practice refer to $(W,f)$ as a noded surface with curves in $\sigma$ {\em pinched}.

\p{Stratification.} We think of $\tsb$ as a stratified space, with $\fkt$ denoting the stratum consisting of surfaces with $k$ curves pinched.  Also, $\ft$ denotes $\tsb-\ts$.  Points in the stratum corresponding to maximally noded surfaces are naturally associated with the pants decompositions that are pinched.  These points are isolated from each other, since there is a unique hyperbolic structure on $S_{0,3}$.

As Teichm\"uller spaces themselves, each connected component $\os$ of $\fkt$ comes equipped with its own metric $d_{\os}$.  On the other hand, $\os$ inherits a metric $d_{\ts}|_{\os}$ as the completion of $\ts$.  Masur proved that the metric tensor on $\ts$ extends continuously to the intrinsic Weil--Petersson metric tensor on strata in $\ft$ \cite{hm}, and Wolpert proved that length minimizing paths in $\fkt$ do not enter $\ts$ \cite{sw2}, and together these facts give the following (see \cite[Lem. 1.3.1]{mw}):

\begin{thm}\label{same}On any connected component $\os$ of $\fkt$, the metrics $d_{\os}$ and $d_{\ts}|_{\os}$ are the same.\end{thm}

\p{Negative curvature.} By work of Tromba and Wolpert the Weil--Petersson metric has negative sectional curvatures and is geodesically convex \cite{ajt} \cite{sw2}.  It then follows from general principles that $\tsb$ is a $\cat$ metric space (see \cite{bh}).

\p{Visual sphere.} Each point $X$ in $\ts$ has a {\em visual sphere} $\vxs$, which is the unit tangent space at $X$.  Due to the non-completeness of the Weil--Petersson metric, certain directions in the visual sphere correspond to finite-length geodesics emanating from $X$ that leave every compact subset of $\ts$.  These {\em finite rays} terminate at noded Riemann surfaces in $\ft$.

As a consequence of the $\cat$ property for $\tsb$, the first author
proved the following density theorem \cite{jb}. 

\begin{thm}\label{jeff}For any $S$ and any $X \in \ts$, the finite rays are dense in the visual sphere $\vxs$.\end{thm} 

The idea is that every point $X$ of $\ts$ is within a uniformly bounded distance of a maximally noded surface $N(X)$ (see \cite{jbpants}).  Thus, given a sequence of points $X_n$ diverging from $X$ along an infinite ray from $X$, a sequence points of $\vxs$ corresponding to the direction determined by the geodesic joining $X$ to $N(X_n)$ converge to the given infinite ray.  Wolpert proved that Theorem~\ref{jeff} has the following corollary \cite{sw}, which immediately implies that a Weil--Petersson isometry is determined by its action on the maximally noded surfaces.

\begin{cor}\label{wolpert}$\tsb$ is the closed convex hull of its
maximally noded surfaces.\end{cor}

Wolpert has in fact already used this corollary to give a simplification of the proof of Masur and Wolf (but not our extension) \cite{sw}.


\section{Edge strata}\label{base}

In this section, we study the stratum of $\tsb$ where all but one of
the curves of a pants decomposition are pinched; in Lemma~\ref{pants},
this stratum is where we will be able to see the edges of $\cps$.  By
Theorem~\ref{same}, each connected component of this stratum is isometric to $\tf$ for $\sooszfo$ equal to either $\soo$ or $\szf$.

In both cases, we have that $\tf$ is topologically $\h$, and $\modf/Z \cong \pglz$, where $Z$ is the center of $\modf$ ($Z = \mathbb{Z}_2$ if $R=\soo$ and $Z = \mathbb{Z}_2 \oplus \mathbb{Z}_2$ if $R=\szf$).  The action of $\modf$ on $\tf$ then coincides with the usual action of $\pglz$ on $\h$ by M\"obius transformations and complex conjugation.

Also, $\ft$ is the set of rational points of $\partial \h$, and $\tfb$ has the horoball topology.  That $\ft$ is discrete agrees with the fact that any noded surface is maximally noded (pants decompositions consist of exactly one curve).  In this way, curves in $\sooszfo$ are naturally identified with the rational numbers, and two curves $p/q$ and $r/s$ (in reduced form) are connected by an edge in $\cpf$ exactly when $|ps-qr|=1$.

\subsection{Farey graph.} Let $\fg$ denote the usual embedding of the \emph{Farey graph} in $\h$:  in the upper half-plane model of $\h$, its vertices are at rational points (including $\infty$), and its edges are hyperbolic geodesics connecting (reduced) rational numbers $p/q$ and $r/s$ whenever $|ps-qr|=1$ (see Figure~\ref{fareypic}).  We will now see that each edge of $\fg$ is a Weil--Petersson geodesic.

The imaginary axis is a geodesic since it is the fixed set of the isometry given by a mapping class corresponding to: \[ \tbt{1}{0}{0}{-1} \] which acts on $\h$ (in the upper half-plane model) by $z \mapsto - \bar{z}$.  The $\pglz$-orbit of this geodesic is $\fg$.  Thus, $\fg$ gives a geodesic between any noded surfaces whose corresponding pants decompositions have distance $1$ in $\cpf$.

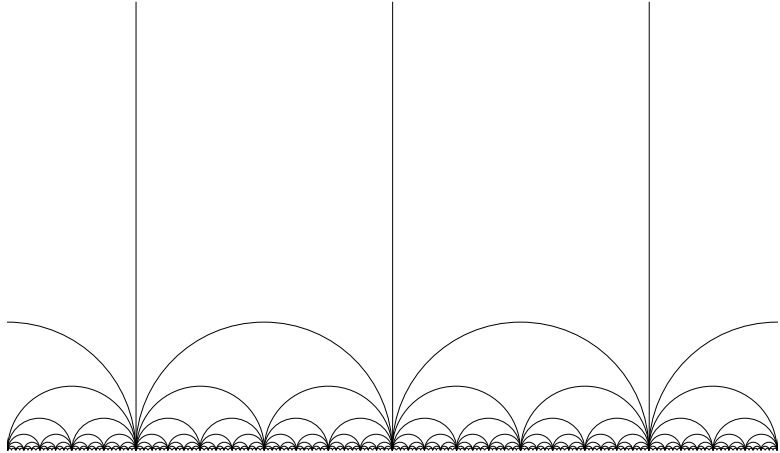
\pics{farey}{The Farey graph $\fg$ in the upper half-plane.}{6}

\subsection{Weil--Petersson vs. pants distance.} We will now show that
the Weil--Petersson metric encodes adjacency in the pants graph $\cpf$.

In the following lemma, we denote a (geodesic) triangle, and any of its edges, by a set of vertices.  A triangle is called a {\em tripod} if any edge is contained in the union of the other two edges.

\begin{lem}\label{ti}Let $(X,d)$ be a $\cat$ metric space, and let $ABC$ be an equilateral triangle in $X$ which is not a tripod, and which has side length $L$.  If $D$ is any point on $BC$, then $d(A,D) > L/2$.\end{lem}

\bpf

Assume that, say, $d(B,D) \leq L/2$.  The triangle inequality implies that $d(A,D) \geq L/2$ and that if $d(A,D) = L/2$, then $d(B,D)=d(C,D)=L/2$.  In the latter case, the uniqueness of geodesics in $\cat$ spaces implies that $ADC=AC$ and $ADB=AB$; thus, $ABC$ is a tripod.

\epf

For the next lemma, $L$ is the Weil--Petersson length of edges of $\fg$, $\pd$ denotes combinatorial distance in $\cpf$, and $\dwp$ denotes Weil--Petersson distance.  

\begin{lem}\label{edges}If $P$ and $P'$ are pants decompositions with $\pd(P,P') > 1$, then the corresponding noded surfaces $W$ and $W'$ have $\dwp(W,W') > L$.\end{lem}

\bpf

Let $g$ be the geodesic from $W$ to $W'$.  Since $\pd(P,P') > 1$, $g$ passes through at least two triangles of $\fg$, and this gives a natural division of $g$ into segments.  By Lemma~\ref{ti}, the segments containing $W$ and $W'$ each have length greater than $L/2$, so $\dwp(W,W') > L$.

\epf

\p{Remark.} $L$ has different values, say $L_{\soo}$ and $L_{\szf}$, depending on the surface; in fact, $L_{\szf} = 2 L_{\soo}$ and more generally the metrics for $\teich(\szf)$ and $\teich(\soo)$ differ by the same factor.  To see this, note that for $R$ either $\soo$ or $\szf$, the quotient of $R$ by $Z(\mcg(R))$ is a genus zero orbifold with one cusp and three cone points of order 2.  Since $\pi:R \to R/Z(\mcg(R))$ is a local isometry, and $\pi$ gives a canonical bijection between $\teich(R)$ and $\teich(R/Z(\mcg(R)))$, it follows from the definition that the metric on $\teich(R)$ is the metric on $\teich(R/Z(R))$ multiplied by the degree of $\pi$.  This fact, however, is not needed for the proof.

\section{Proof of Main Theorem}

Let $S$ be any surface, other than $S_{0,3}$, with negative Euler characteristic, let $I \in \isom(\ts)$, and let $\ib$ be the natural extension of $I$ to the completion $\tsb$.  We will use $\fgek$ to denote $\displaystyle \bigcup_{i \geq k} \fri$.  

The following lemma is also proven by Masur and Wolf \cite[Lem. 1.3.4, Lem. 1.3.5]{mw}; however, the proof given here is more elementary, relying only on topology and not at all on curvature.

\begin{lem}\label{strata}$\ib(\ft) \subset \ft$, $I$ and $\ib$ are surjective, and $\ib$ preserves strata.\end{lem}

\bpf

For the first statement, suppose that $\ib(W) \in \ts$.  Then $\ib(W)$ has a compact neighborhood $K$.  Let $\kb = K \cap \ib(\tsb)$.  We have that $\kb$ is compact, since $\ib(\tsb)$ is closed (it is an isometric embedding of a complete space).  But this implies that $W$ has a compact neighborhood (namely $\ib^{-1}(\kb)$), and so $W \in \ts$.

For surjectivity of $I$, it suffices to show that $I$ is proper (any proper embedding of a manifold without boundary into itself is surjective).  Indeed, let $K$ be a compact subset of $\ts$.  As above, $\kb = K \cap \ib(\tsb)$ is compact.  Since $\ib(\ft) \subset \ft$, $\kb = K \cap I(\ts)$.  Hence, $I^{-1}(\kb)$ is compact.  Thus, $I$ is surjective, and it immediately follows that $\ib$ is surjective.

For the last statement, we inductively show that $W \in \fkt$ if and only if $\ib(W) \in \fkt$.  The base case is $\fzt = \ts$:  by definition, $\ib(W) \in \ts$ whenever $W \in \ts$, and we have just shown that $\ib(\ft) \subset \ft$.

Now suppose $\ib$ preserves $\fkmo$.  Then $\ib$ restricts to an isometry of $\fgek$, which is a homeomorphism since $\ib$ is surjective.  The points of $\fkt$ are characterized by having compact neighborhoods in $\fgek$.  Since $\ib|_{\fgek}$ is a homeomorphism, $\fkt$ is preserved.

\epf

We now combine Lemmas~\ref{edges} and~\ref{strata} to bridge the gap between isometries of $\ts$ and automorphisms of $\cps$.

\begin{lem}\label{pants}$\ib$ induces $\ibs \in \aut(\cps)$.\end{lem}

\bpf

By Lemma~\ref{strata}, we have that $\ib$ induces $\ibs \in \aut(\cpos)$.  It remains to show that $\ibs$ preserves edges in $\cps$.

Points in the highest level stratum $\fmk$ correspond to pants
decompositions.  If $\os$ is any connected component of $\fmok$, then
$\os$ and $\ib(\os)$ are both copies of either $\tsoo$ or $\tszf$.  As
with any connected component of $\fmok$, $\os$ comes equipped with a
Farey graph $\fg_{\os}$.  The lengths of the edges of $\fg_{\os}$ and
$\fg_{\ib(\os)}$ must be the same, since the smallest length of a
geodesic connecting noded surfaces is an isometry invariant (in fact,
by the remark at the end of Section~\ref{base}, the surfaces
corresponding to $\os$ and $\ib(\os)$ are the same).  By
Lemma~\ref{edges}, Lemma~\ref{strata} and Theorem~\ref{same},
$\ib(\fg_{\os}) = \fg_{\ib(\os)}$; but this is exactly saying that
$\ibs \in \aut(\cps)$.

\epf

We are now ready to prove the Main Theorem.

\bpf

Let $I \in \isom(\ts)$.  By Lemma~\ref{pants}, $I$ induces $\ibs \in \aut(\cps)$.  By Theorem~\ref{dan}, $\ibs$ is induced by a mapping class $f$.  By Corollary~\ref{wolpert}, $f$ induces $I$.

Also, by Corollary~\ref{wolpert}, the kernel of the map $\eta: \mods \to \isom(\ts)$ is the same as the kernel of $\theta: \mods \to \aut(\cps)$.

\epf

\bibliographystyle{plain}
\bibliography{wp}

\end{document}